\documentclass[12pt, reqno]{amsart}
\usepackage{amsfonts}
\usepackage{amsmath}
\usepackage{amssymb}
\usepackage{eufrak}
\usepackage{mathrsfs}
\usepackage{bbding}

\numberwithin{equation}{section} 

\input xy
\xyoption{all}
\title{Brill-Noether general curves on Knutsen K3 surfaces}

\author{Maxim Arap}  \author{Nicholas Marshburn}
\address{Maxim Arap \\ Johns Hopkins University\\ Department of Mathematics \\  404 Krieger Hall \\ 3400 N. Charles Street \\ Baltimore, MD 21218\\ USA}
\email{marap@math.jhu.edu}
\address{Nicholas Marshburn \\ Johns Hopkins University\\ Department of Mathematics \\  404 Krieger Hall \\ 3400 N. Charles Street \\ Baltimore, MD 21218\\ USA}
\email{marshbur@math.jhu.edu}

\begin{document}

\maketitle

\newcommand{\bP}{\mathbb{P}}
\newcommand{\bC}{\mathbb{C}}
\newcommand{\mc}{\mathcal}
\newcommand{\ra}{\rightarrow}
\newcommand{\thra}{\twoheadrightarrow}
\newcommand{\mb}{\mathbb}
\newcommand{\mrm}{\mathrm}
\newcommand{\p}{\prime}
\newcommand{\ms}{\mathscr}
\newcommand{\pl}{\partial}
\newcommand{\ti}{\tilde}
\newcommand{\wti}{\widetilde}
\newcommand{\ol}{\overline}
\newcommand{\ul}{\underline}
\newcommand{\Sp}{\mrm{Spec}\,}

\newtheorem{pro1}{Proposition}[section]
\newtheorem{pro2}[pro1]{Proposition}
\newtheorem{pro3}[pro1]{Proposition}
\newtheorem{pro4}[pro1]{Proposition}
\newtheorem{pro5}[pro1]{Proposition}
\newtheorem{pro6}[pro1]{Proposition}
\newtheorem{pro7}[pro1]{Proposition}

\newtheorem{thm1}[pro1]{Theorem}
\newtheorem{thm2}[pro1]{Theorem}
\newtheorem{thm3}[pro1]{Theorem}
\newtheorem{thm4}[pro1]{Theorem}
\newtheorem{thm5}[pro1]{Theorem}

\newtheorem{cor1}[pro1]{Corollary}
\newtheorem{cor2}[pro1]{Corollary}
\newtheorem{cor3}[pro1]{Corollary}
\newtheorem{cor5}[pro1]{Corollary}

\newtheorem{lem1}[pro1]{Lemma}
\newtheorem{lem2}[pro1]{Lemma}
\newtheorem{lem3}[pro1]{Lemma}
\newtheorem{lem4}[pro1]{Lemma}
\newtheorem{lem5}[pro1]{Lemma}
\newtheorem{lem6}[pro1]{Lemma}
\newtheorem{lem7}[pro1]{Lemma}
\newtheorem{lem8}[pro1]{Lemma}
\newtheorem{lem9}[pro1]{Lemma}

\newtheorem{rem1}[pro1]{Remark}
\newtheorem{rem2}[pro1]{Remark}
\newtheorem{rem3}[pro1]{Remark}

\begin{abstract} This article classifies Knutsen K3 surfaces all of whose hyperplane sections are irreducible and reduced. As an application, this gives infinite families of K3 surfaces of Picard number two whose general hyperplane sections are Brill-Noether general curves. 
\end{abstract} 

\thispagestyle{empty}
\section{Introduction}

In \cite[Thm.1.1]{Knu}, one may find a classification of triples of integers $(n,d,g)$ such that there exists a smooth K3 surface of degree $2n$ in $\bP^{n+1}$ which contains a smooth curve of degree $d$ and genus $g$. Moreover, for each such triple Knutsen constructed a K3 surface, which we shall denote by $S_{n,d,g}$ and call \emph{Knutsen K3 surface}, containing a smooth curve of degree $d$ and genus $g$. By construction, Knutsen K3 surfaces have Picard number either one or two.  

A smooth irreducible projective curve $C$ is said to be \emph{Brill-Noether general} if the Petri map
$$H^0(L)\otimes H^0(\omega_C\otimes L^*) \ra H^0(\omega_C)$$
defined by multiplication is injective for every line bundle $L$ on $C$. A theorem of Lazarsfeld's says that if $S$ is a smooth K3 surface all of whose hyperplane sections are irreducible and reduced then a general hyperplane section is a Brill-Noether general curve. In particular, general sections of K3 surfaces of Picard number one are Brill-Noether general.   

It is well known that generic K3 surfaces of a given degree have Picard number one and Lazarsfeld's theorem immediately applies. However, in practice one sometimes has to deal with non-generic K3 surfaces of Picard number $\ge 2$ and the question of whether such surfaces have Brill-Noether general sections arises. For instance, in \cite{ACM} the question of whether a given (Knutsen) K3 surface $S$ embeds in a certain Fano threefold is reduced to the question of whether $S$ has a Brill-Noether general section. 

In this article we determine Knutsen K3 surfaces $S_{n,d,g} \subset \bP^{n+1}$ of Picard number two all of whose hyperplane sections are irreducible and reduced (Theorem \ref{thm1}). As a corollary, we obtain a numerical condition on $n,d,g$ that guarantees that $S_{n,d,g}$ has a  Brill-Noether general hyperplane section (Corollary \ref{cor1}).  

\vspace{1cm}

\textbf{Notation and conventions.} We work over the field $\mb{C}$ of complex numbers. By a curve we shall mean a reduced scheme over $\mb{C}$ of dimension one. All K3 surfaces are assumed to be smooth and projective. For a real number $r$, the symbol $\lceil r \rceil$ denotes the smallest integer $\ge r$. The symbol $\sim$ denotes linear equivalence of divisors. 

\section{Irreducible and reduced hyperplane sections}
In what follows we assume $n\ge 2$. By \cite{Knu}, Knutsen K3 surface $S_{n,d,g} \subset \bP^{n+1}$ of Picard number two with hyperplane section $H$ and smooth curve $C \subset S_{n,d,g}$ of degree $d$ and genus $g$ has $\mrm{Pic}(S_{n,d,g})= \mb{Z}H\oplus \mb{Z}C$ with the following intersection matrix
$$\begin{bmatrix}
H^2 & H\cdot C\\
C\cdot H & C^2
\end{bmatrix} = 
\begin{bmatrix} 
2n & d \\
d & 2g-2
\end{bmatrix}.
$$
The following proposition reduces the main question of this article to a system of diophantine inequalities.  
 
\begin{pro1} \label{pro1}
Let $S:=S_{n,d,g} \subset \bP^{n+1}$ be a Knutsen K3 surface with $\mrm{Pic}(S)=\mb{Z}H\oplus \mb{Z}C$. The following conditions are equivalent:
\begin{enumerate}
\item The linear system $|H|$ contains a reducible or non-reduced member. 
\item There exists an irreducible curve of degree $\le n$ on $S$. 
\item There exist integers $a,b$ satisfying 
\begin{equation*}
\left\{ 
\begin{array}{lr}
0 < 2na+bd \le n  & \mathrm{        (I)}\text{ }\\
na^2+dab+(g-1)b^2\ge -1 & \mathrm{        (II)}.
\end{array} 
\right.
\end{equation*} 
\end{enumerate}
\end{pro1}
\begin{proof} $(1)\!\!\implies\!\! (2)\colon$ Since $H^2=2n$ then for any splitting $H\sim D_1+D_2$ we must have $\deg D_1 \le n$ or $\deg D_2 \le n$. Say $\deg D_1 \le n$, then an irreducible and reduced component of $D_1$ gives an irreducible curve of degree $\le n$ on $S$. 

$(2)\!\! \implies \!\!(3)\colon$ Say $D$ is an irreducible curve of degree $\le n$ on $S$ and write $D \sim aH+bC$ for some $a,b \in \mb{Z}$. The inequalities (I) and (II) follow immediately from $0 < \deg D \le n$ and $D^2 \ge -2$, respectively. 

$(3)\!\! \implies \!\!(1)\colon$ Let $a,b$ be integers satisfying (I) and (II), and let $D := aH+bC$. The inequalities (I) and (II) imply $0 < \deg D \le n$ and $D^2 \ge -2$. Using Hirzebruch-Riemann-Roch formula and $D^2 \ge -2$ we may check that $|D|$ is non-empty. Let $E$ be an irreducible and reduced component of an effective divisor in $|D|$. We have $\deg E \le n$, and therefore, $E$ is contained in a hyperplane $\bP^n \subset \bP^{n+1}$. Let $E^\prime$ be the complement of $E$ in the intersection $H:=\bP^n \cap S$ then we have a splitting $H=E+E^\p$. 
\end{proof}

\begin{lem1} \label{lem1}
Let $S:=S_{n,d,g}$ be a Knutsen K3 surface with $\mrm{Pic}(S)=\mb{Z}H\oplus \mb{Z}C$. Suppose $b\in \mb{Z}$ and $a:=\big\lceil{-}\frac{bd}{2n}\big\rceil$ are such that $D:=aH+bC$ satisfies $0 < \deg D \le n$ and $D^2 \ge -2$. Let $r$ be the residue of $d$ modulo $2n$. Then one of the following two conditions holds:
\begin{enumerate}
\item $r \le n$, $(a_1H+C)^2\ge -2$ with $a_1=\big\lceil{-}\frac{d}{2n}\big\rceil$
\item $r > n$, $(a_1H-C)^2\ge -2$  with $a_1=\big\lceil\frac{d}{2n}\big\rceil$.
\end{enumerate}
\end{lem1}
\begin{proof} Let $\delta = r$ if $r\le n$ and $\delta=2n-r$ if $r>n$. Also, let $\epsilon$ be the residue of $bd$ modulo $2n$. We may check that $\deg D \le n$ implies that $0 \le \epsilon \le n$.

First, let us consider the case $r \le n$. We may check that the inequality $D^2\ge -2$ is equivalent to the condition
$$\frac{\epsilon^2-d^2b^2}{4n}+(g-1)b^2 \ge -1.$$
Since $\delta b\equiv db \equiv \epsilon \text{ }(\text{mod } 2n)$ and $0 \le \epsilon \le n$ then $|\delta b| \ge \epsilon$, and therefore, we have 
$$-1 \le \frac{\epsilon^2-d^2b^2}{4n}+(g-1)b^2 \le \frac{\delta^2b^2-d^2b^2}{4n}+(g-1)b^2.$$
Since $\deg D>0$ then $b\ne 0$. Dividing by $b^2$ we obtain 
$$-1 \le -\frac{1}{b^2}\le \frac{\delta^2-d^2}{4n}+g-1=\frac{1}{2}(a_1H+C)^2,$$
where $a_1=\frac{\delta-d}{2n}=\big\lceil{-}\frac{d}{2n}\big\rceil$, which gives the desired inequality $(a_1H+C)^2 \ge {-}2$. 

In the case when $r>n$ an analogous calculation with $\delta=2n-r$ shows that $(a_1H-C)^2 \ge -2$ with $a_1=\frac{\delta-d}{2n}=\big\lceil\frac{d}{2n}\big\rceil$.
\end{proof}

\begin{thm1} \label{thm1}
Knutsen K3 surface $S_{n,d,g} \subset \bP^{n+1}$ of Picard number two has a reducible or non-reduced hyperplane section if and only if $g\ge \frac{d^2-\delta^2}{4n}$, where $\delta$ is the distance from $d$ to the nearest integer multiple of $2n$. 
\end{thm1}
\begin{proof} By Proposition \ref{pro1}, $S_{n,d,g} \subset \bP^{n+1}$ has a reducible or non-reduced hyperplane section if and only if there exist integers $a,b$ satisfying inequalities (I) and (II). We may check that if $a,b\in \mb{Z}$ give a solution to (I) then $a=\big\lceil{-}\frac{bd}{2n}\big\rceil$. Therefore, by Lemma \ref{lem1}, the integers $a,b$ satisfying (I) and (II) exist if and only if one of the conditions $(1)$ or $(2)$ of Lemma \ref{lem1} holds. It is easily seen that each of the conditions $(1)$ and $(2)$  of Lemma \ref{lem1} is equivalent to the requirement $g\ge \frac{d^2-\delta^2}{4n}$, where $\delta$ is the distance from $d$ to the nearest integer multiple of $2n$. 
\end{proof}

\begin{cor1} \label{cor1}
Let $S_{n,d,g} \subset \bP^{n+1}$ be a Knutsen K3 surface of Picard number two and let $\delta$ be the distance from $d$ to the nearest integer multiple of $2n$. If $g < \frac{d^2-\delta^2}{4n}$ then a general hyperplane section of $S_{n,d,g}$ is a Brill-Noether general curve. 
\end{cor1}
\begin{proof}
The proof follows immediately from Theorem \ref{thm1} and Lazarsfeld's theorem. 
\end{proof}
%\section*{Acknowledgements} 


\begin{thebibliography}{99}
\bibitem[ACM]{ACM} Arap, M., Cutrone, J., Marshburn N.: \emph{On the existence of certain weak Fano threefolds of Picard number two.} arXiv:1112.2611v1.
\bibitem[Knu]{Knu} Knutsen, A.L.: \emph{Smooth curves on projective $K3$ surfaces.} Math. Scand. \textbf{90} (2002) 215--231.
\bibitem[Laz]{Laz} Lazarsfeld, R.: \emph{Brill-Noether-Petri without degenerations.} J.
Diff. Geom. \textbf{23} (1986), no. 3, 299--307.
\end{thebibliography}
\end{document}